\documentclass{amsart}
\usepackage{amssymb,latexsym,graphics,amsfonts}

\swapnumbers
 \theoremstyle{plain}

 \theoremstyle{definition}
 
 \theoremstyle{remark}

 \newcommand{\set}[1]{\left\{#1\right\}}

 \newcommand{\cal}[1]{\mathcal{#1}}

\begin{document}

\title[First cohomology on weighted semigroup algebras] {First
cohomology on weighted semigroup algebras}

\author{M. Eshaghi Gordji}
\address{Department of Mathematics,
University of Semnan, Semnan, Iran} \email{maj\_ess@Yahoo.com \&
meshaghi@semnan.ac.ir}

\author{ F. Habibian}
\address{Department of Mathematics,
Isfahan University, Isfahan, Iran}
\email{fhabibian@math.ui.ac.ir}

\author{A. Rejali}
\address{Department of Mathematics,
Isfahan University, Isfahan, Iran}\email{rejali@math.ui.ac.ir}



\dedicatory{}


\smallskip

\begin{abstract}
The aim of this work is to generalize Johnson's techniques in order
to apply them to establish a bijective correspondence between
$S$-derivations and continuous derivations on $M_a(S,\omega),$ where
$S$ is a locally compact foundation semigroup with identity $e$, and
$\omega$ is a weight function on $S$, and apply it to find a
necessary condition for amenability of weighted group algebras.
\end{abstract}
\maketitle


\section{Introduction}
Let $G$ be a locally compact group and $L^1(G)$ be the group algebra
of $G$. We shall regard an neo-unital Banach $L^1(G)$-module as a
unital Banach $M(G)$-module, and hence as a $l^1(G)$-module, in this
way.
 B. E. Johnson has established a bijective
correspondence between continuous derivations
$D:L^1(G)\longrightarrow E^\prime$ and $G$-derivations
$\varrho:G\longrightarrow E^\prime$ and further, he proved that
inner derivations correspond to inner $G$-derivations, and
$H^1(L^1(G),E^\prime)=\{0\}$ if and only if every $G$-derivations
into $E^\prime$ is inner where $E$ is a neo-unital Banach
$L^1(G)$-module (see [11] and [12]).
Throughout, $S$ denotes a locally compact Hausdorff topological
semigroup. We shall assume that there is a $weight$ $function$
 $\omega$ on $S$ .
 By this mean that $\omega:S\longrightarrow
[0,\infty)$ with the properties that $\omega$ is continuous
 and that both $\omega$ and $1/\omega$ are bounded on compact subsets
of $S$ and that $\omega(xy)\leq\omega(x)\omega(y)$ for all $x,y\in
S$. A complex-valued function $f$ on $S$ is called $\omega$-bounded
if there is $k>0$ such that $\mid f(x)\mid\leq k\omega(x)$ for all
$x\in S$. Recall that $M_a(S)$ denotes the space of all measures
$\mu\in M(S)$, the space of all bounded complex Radon measures on
$S$, for which the mappings $x\longrightarrow \delta_x*\mid \mu
\mid$ and $x\longrightarrow \mid \mu \mid*\delta_x$ from $S$ into
$M(S)$ are weakly continuous, where $\delta_x$ denotes the Dirac
measure at $x$(see [5], or $\tilde{L}(S)$ as in [1]). A topological
semigroup $S$ is called $foundation$ $ semigroup$ if $S$ coincides
with the closure of $\bigcup\{supp(\mu):\mu\in M_a(S)\}$. The space
of all complex Radon measures $\mu$ on $S$ for which $\omega\mu\in
M(S)$ is denoted by $M(S,\omega)$. With the convolution
$(\mu*\nu)(f)=\int_S\int_S f(xy)d\mu(x)d\nu(y),$ for each continuous
complex-valued function $f$ on $S$ with compact support (so for
every $\omega$-bounded Borel measurable function $f$ on $S$), and
the norm $\parallel \mu
\parallel_\omega=\parallel \omega\mu
\parallel=\int_S\omega(t)d\mid \mu \mid(t)$, $M(S,\omega)$ defines
a convolution measure algebra. Also, the space of all measures
$\mu\in M(S,\omega)$ such that $\omega\mu\in M_a(S)$ is denoted by
$M_a(S,\omega)$ (see [9], or $\tilde{L}(S,\omega)$ as in [10]). By
theorem 4.6 of [10], $M_a(S,\omega)$ is a two-sided L-ideal of
$M(S,\omega)$. We denote by $L^\infty(S,;M_a(S))$ the set of all
bounded complex-valued $\mu$-measurable ($\mu\in M_a(S)$) functions
on $S$, and we also define
$L^\infty(S;M_a(S,\omega))=\{f:\frac{f}{\omega}\in
L^\infty(S;M_a(S))\}.$ We recall from lemma 2 of [8] that, if $S$ is
a foundation semigroup with identity and with a weight function
$\omega$, then $L^\infty(S;M_a(S,\omega))$ with the norm $\parallel
.
\parallel_\omega$ where $\parallel f
\parallel_{\omega}=\parallel f/\omega\parallel_\infty$ for $f\in
L^\infty(S;M_a(S,\omega))$ and the multiplication $(f\odot
g)(x)=\frac{f(x)g(x)}{\omega(x)}$ for $f,g\in
L^\infty(S;M_a(S,\omega)),$ and the involution $f^*=\bar{f}$,
defines a commutative $C^*$-algebra. Moreover, the mapping
$f\longmapsto \tau_f :L^\infty(S;M_a(S,\omega))\longrightarrow
M_a(S,\omega)$, defined by \[\tau_f(\mu)=\int_Sf(x)d\mu(x)\qquad
(\mu\in M_a(S,\omega)),\] is an isometric isomorphism from
$L^\infty(S;M_a(S,\omega))$ onto the dual space of $M_a(S,\omega)$.
\\Let $\cal A$ be a Banach algebra and $X$ be a Banach $\cal
A$-module, then $X^*$ is a Banach $\cal A$-module if for each $a\in
\cal A$ and let $x\in X$ and $x^*\in X^*$ we define
\[\langle x,ax^*\rangle=\langle xa,x^*\rangle, \hspace {1 cm} \langle
x,x^*a\rangle=\langle ax,x^*\rangle.\]
 Let $X$ be a Banach $\cal A$-module, then a derivation from
$\cal A$ into $X$ is a continuous linear operator $D$, if for every
$a,b\in \cal A,$
\[D(ab)=D(a).b+a.D(b).\] If $x\in X$ and we define
$D$ by \[D(a)=a.x-x.a\qquad (a\in \cal A),\] then $D$ is a
derivation. Derivations of this form are called inner derivations.
 A Banach algebra $\cal A$  is amenable if every derivation from $\cal
A$ into every dual $\cal A$-module is inner, equivalently if
$H^1(\cal A,X^*)=\{0\}$ for every Banach $A$-module $X$, where
$H^1(\cal A,X^*)$ is the first cohomology group from $\cal A$ with
coefficients in $X^*$[11] (see [3] and [14], for more
details).\\
\section{weighted measure algebras}
 Let $S$ be a locally compact foundation semigroup with
identity $e$, and  let $\omega$ be a weight function on $S$ such
that $\omega(e)=1$. Choose a fixed neighbourhood $U_0$ of $e$. Let
$\Lambda$ denote the collection of all compact neighbourhoods of $e$
contained in $U_0$, ordered by inclusion. For every $\lambda\in
\Lambda$ let $\mu_\lambda$ be a measure in $M_a(S)$ such that
$\mu_\lambda\geq0$, $\parallel\mu_\lambda\parallel=1$, and its
support contained in $\lambda$. Put
$\nu_\lambda=\frac{\mu_\nu}{\omega(x)}$. Then $\parallel
\nu_\lambda\parallel_\omega=1$. Let $\mu\in M_a(S,\omega)$. Then
from the $\parallel.\parallel_\omega$-continuity of the mapping:
$x\longrightarrow \mu*\delta_x$ from $S$ into $M_a(S,\omega)$ and
the continuity of $\omega$ at $e$ it follows that for every
$\varepsilon>0$ there exists a $\lambda_1\in \Lambda$ such that
$\parallel\mu*\delta_x-\mu\parallel_\omega<\varepsilon$ and
$\mid\omega(x)-1\mid<\varepsilon$ for all $x\in \lambda_1$. Thus for
all $\lambda\geq\lambda_1$
\[\parallel\mu*\nu_\lambda-\mu\parallel_\omega<\varepsilon(1+\parallel\mu\parallel_\omega).\]
Then $(\nu_\lambda)$ is a bounded approximate identity for
$M_a(S,\omega)$ such that $\parallel \nu_\lambda\parallel_\omega=1$,
$\nu_\lambda\geq0$, and furthermore its support contained in
$\lambda$ for each $\lambda\in \Lambda$ (see [8]).
\paragraph{\bf Definition 2.1.}Let $S$ be a locally compact foundation
semigroup with identity, and  let $\omega$ be a weight function on
$S$. A left weighted Banach $S$-module $E$ is a Banach
space $E$ which is a left $S$-module such that\\
(1) for all $x$ in $E$ the map $s\longmapsto s.x$ is continuous from
$S$ into
$E$.\\
(2) There is $K>0$ such that $\parallel s.x\parallel\leq
K.\omega(s).\parallel x\parallel$ for all $s\in S$ and $x\in E$.\\
In a similar way we define right weighted Banach $S$-module and
two-sided weighted Banach $S$-modules.
\paragraph{\bf Definition 2.2.} If a Banach algebra $\cal A$ is
contained as a closed ideal in another Banach algebra $\cal B$, then
the $strict$ topology or $strong$ operator topology $(so)$ on $B$
with respect to $\cal A$ is defined through the family of seminorms
$(p_a)_{a\in \cal A}$, where $p_a(b)=\parallel ba
\parallel+\parallel ab \parallel$ for every $ b\in \cal B.$

Let $S$ be a locally compact foundation semigroup with identity, and
let $\omega$ be a weight function on $S$, we say that a net
$(\mu_\alpha)_\alpha$ in $M(S,\omega)$ converges to $\mu\in
M(S,\omega)$ in (so), if $\mu_\alpha*\nu\longrightarrow \mu*\nu
\hspace{0.5 cm}$and $\hspace{0.5 cm} \nu*\mu_\alpha \longrightarrow
\nu*\mu$, for every $\nu\in M_a(S,\omega)$.
\paragraph{\bf Lemma 2.3.}Let $S$ be a locally compact foundation
semigroup with identity,
 and  let $\omega$ be a weight function on $S$. Then
$F\in(M(S,\omega),so)^\prime$ if and only if there exists subsets
$\set{\nu_1,...,\nu_n}$ and $\set{\eta_1,...,\eta_m}$ in
$M_a(S,\omega)$, and $\set{\varphi_1,...,\varphi_n}$ and
$\set{\psi_1,...,\psi_m}$ in
$M_a(S,\omega)^\prime=L^\infty(S;M_a(S,\omega))$ such that for each
$\mu\in M(S,\omega)$;\[F(\mu)=\sum_{i=1}^n
\varphi_i(\nu_i*\mu)+\sum_{j=1}^m \psi_j(\mu*\eta_j),\quad\] where
\[\varphi_i(\nu_i*\mu)=\int\varphi_i d(\nu_i*\mu).\]
\paragraph{\bf Proof.}
The sufficient  condition is evident.
 Let $F\in(M(S,\omega),so)^\prime$
then by theorem 3.3 of [2] there is $\set{\nu_1,...,\nu_n}$ of
$M_a(S,\omega)$ which \[\mid F(\mu)\mid\leq\sum_{i=1}^n (\parallel
\nu_i*\mu
\parallel_\omega+\parallel \mu*\nu_i \parallel_\omega).\] We define
\[\Delta=\set{(\nu_1*\mu,...,\nu_n*\mu,\mu*\nu_1,...,\mu*\nu_n): \mu\in
M(S,\omega)}\] and $\gamma:\Delta\longrightarrow \mathbb{C}$ by
$\gamma(T)=F(\mu)$. Clearly $\Delta\subseteq\bigoplus_{i=1}^n
M_a(S,\omega)$ and $\gamma$ is well-defined and bounded. By
Hahn-Banach theorem there is a bounded functional $\Gamma$ on
$\bigoplus_{i=1}^n M_a(S,\omega)$ such that
$\Gamma\mid_\Delta=\gamma$.\\Now for all $1\leq i\leq n$ and $1\leq
j\leq 2$ we define $\Phi_{ij}\in
M_a(S,\omega)^\prime=L^\infty(S;M_a(S,\omega))$ by
\[\Phi_{ij}(\nu)=\Gamma(0,...,0,\nu_{(j-1)n+i},0,...,0) \quad (\nu\in
M_a(S,\omega)),\]and hence for each $\mu\in M(S,\omega)$ we have
\[F(\mu)=\gamma(\nu_1*\mu,...,\nu_n*\mu,\mu*\nu_1,...,\mu*\nu_n)\]
 \[\qquad=\Gamma(\nu_1*\mu,...,\nu_n*\mu,\mu*\nu_1,...,\mu*\nu_n)\]
 \[\qquad=\sum_{i=1}^n \Phi_{i1}(\nu_i*\mu)+\sum_{i=1}^n
 \Phi_{i2}(\mu*\nu_i).\] and proof is complete.\hfill$\blacksquare~$
\paragraph{\bf Theorem 2.4.} Let $S$ be a locally compact foundation
semigroup with identity, and  let $\omega$ be a weight function on
$S$. Then $l^1(S,\omega)$ is $(so)$ dense in $M(S,\omega)$.
\paragraph{\bf Proof.}
By theorem 4.7 of [13],
$\overline{l^1(S,\omega)}^{w^*}=(^\perp(l^1(S,\omega)))^\perp$,
where \[^\perp(l^1(S,\omega)=\{\varphi\in l^1(S,\omega)^\prime:
\varphi(f)=0, \forall f\in l^1(S,\omega)\},\]and hence
$\overline{l^1(S,\omega)}^{w^*}\cap M(S,\omega)=M(S,\omega)$. Now if
$\mu\in M(S,\omega)$ then there exists a net $(\mu_\alpha)_\alpha$
in $l^1(S,\omega)$ such that $\mu_\alpha\longrightarrow \mu$ in
$w^*$-topology. Let $F\in
M_a(S,\omega)^\prime=L^\infty(S;M_a(S,\omega))$ then for each
$\nu\in M_a(S,\omega)$, we have
\[F(\mu_\alpha*\nu)=\int_G F(x)d(\mu_\alpha*\nu)(x)=
\int_G\int_G F(xy)d\mu_\alpha(x)d\nu(y)\]
\[\qquad=\int_G(\int_G
F(xy)d\nu(y))d\mu_\alpha(x)\longrightarrow\int_G(\int_G
F(xy)d\nu(y))d\mu(x)\]
\[\qquad=\int_G\int_G
F(xy)d\mu(x)d\nu(y)=\int_G F(x)d(\mu*\nu)(x)=F(\mu*\nu)\] then by
lemma 2.3, for any $F\in(M(S,\omega),so)^\prime$  we have
$F(\mu_\alpha)\longrightarrow F(\mu)$. This imply that $\mu$ is in
the closure of $l^1(S,\omega)$ with weak topology. Since
$l^1(S,\omega)$ is convex, by theorem 3.12 of [13],  $\mu$ is in
closure of $l^1(S,\omega)$ in $(so)$ and the proof is complete.
\hfill$\blacksquare~$
\paragraph{\bf Corollary 2.5.}
Let $S$ be a locally compact foundation semigroup with identity, and
let $\omega$ be a weight function on $S$. Let $E$ be a neo-unital
Banach $M_a(S,\omega)$-module, and let $D\in
\emph{Z}^1(M_a(S,\omega),E^\prime)$. Then $E$ is a neo-unital Banach
$M(S,\omega)$-module, and there is a unique $\tilde{D}\in
\emph{Z}^1(M(S,\omega),E^\prime)$ that extends $D$ and is continuous
with respect to so-topology on $M(S,\omega)$ and the $w^*$-topology
on $E^\prime$. In particular, $\tilde{D}$ is uniquely determined by
its values on $\{\delta_s:s\in S\}$.
\paragraph{\bf Proof.} The proof of the first part follows from 2.1.6
of [14]. By 2.4, $l^1(S,\omega)$ is strictly dense in $M(S,\omega)$,
which shows that $\tilde{D}$ is uniquely determined through its
values on $l^1(S,\omega)$ and hence on $\{\delta_s:s\in S\}$.
\hfill$\blacksquare~$
\paragraph{\bf Definition 2.6.}Let $S$ be a locally compact foundation
semigroup with identity, and  let $\omega$ be a weight function on
$S$. Let $E$ be a weighted $S$-module. A map
$\varrho:S\longrightarrow E$ is a weighted $S$-derivation if,\\
(1) $\varrho(st)=\varrho(s).t+s.\varrho(t)$, for each $s,t\in
S$.\\
(2) $\varrho$ is $\omega$-bounded.\\ Furthermore, $\varrho$ is inner
if there exists $x\in E$ such that $\varrho(s)=s.x-x.s$, for all
$s\in S$.
\paragraph{\bf Lemma 2.7.}Let $S$ be a locally compact foundation
semigroup with identity and  let $\omega$ be a weight function on
$S$. $X$ is a weighted $S$-module if and only if, $X$ is a
neo-unital Banach $M_a(S,\omega)$-module.
\paragraph{\bf Proof.} Let $\mu\in M(S,\omega)$. We define
\[\mu.x=\int_S s.x
d\mu(s),\qquad  and \qquad x.\mu=\int_S x.s d\mu(s).\] We only check
that $X$ is a neo-unital left Banach $M_a(S,\omega)$-module by
$\mu.x$ as an action. Since for some $K>0$
\[\int_S\parallel s.x\parallel d\mid \mu \mid(s)\leq K.\parallel
x\parallel\int_S\omega(s)d\mid \mu \mid(s)=K.\parallel
x\parallel.\parallel \mu \parallel_\omega<\infty,\] then $\mu.x$
 is well-defined. We have
\[\mu.(\nu.x)=\int_S s\int_S t.x
d\nu(t)d\mu(s)=\int_S\int_S(s.t).x d(\mu\otimes\nu)(s,t)=\int_S u.x
d\mu\nu(u)=(\mu.\nu).x,\] and \[\parallel \mu.x\parallel \leq
\int_S\parallel s.x\parallel d\mid \mu \mid(s)\leq K.\parallel
x\parallel\int_S\omega(s)d\mid \mu \mid(s)=K.\parallel
x\parallel.\parallel \mu \parallel_\omega.\]It is clear that the map
$\mu \longmapsto \mu.x$ is continuous from $M(S,\omega)$ into $X$.
Then $X$ is an $M(S,\omega)$-module, and an $l^1(S,\omega)$-module
by restriction of these operations to the absolutely continuous
measures. Let $(\nu_\lambda)$ be the bounded approximate identity
for $M_a(S,\omega)$ described at the first of this section. Let
$\varepsilon>0$ be given. Since the map $s\mapsto s.x$ is
continuous, there is a neighbourhood $U_\varepsilon$ such that
\[s\in U_\varepsilon\Longrightarrow \parallel s-s.x \parallel<
\varepsilon.\] Then for all $\lambda \geq \lambda_\varepsilon$ we
have $supp(\nu_\lambda)\subseteq U_\varepsilon$. Now let $\lambda
\geq \lambda_\varepsilon$, then \[\parallel\int_S
s.xd\nu_\lambda(s)-x\parallel\leq \int_S \parallel s.x-x
\parallel d\nu_\lambda(s)=\int_{U_\varepsilon} \parallel s.x-s
\parallel d\nu_\lambda(s)\leq \varepsilon
\int_{U_\varepsilon}d\nu_\lambda(s)=\varepsilon,\] for all $x\in X$
, so $x\in \overline{M_a(S,\omega)X},$  and then
$X=\overline{M_a(S,\omega)X}$. But $M_a(S,\omega)$ has a bounded
approximate identity, then by Cohen factorization theorem
$X=M_a(S,\omega)X$ and hence $X$ is neo-unital.\\Conversely, let $X$
be a neo-unital Banach $M_a(S,\omega)$-module. Then by 2.5, $X$ is a
neo-unital $M(S,\omega)$-module. We define $s.x=\delta_s.x$ for
$s\in S, x\in X$. We claim that, $X$ is a weighted Banach $G$-module
by this action. It is enough to show that there is $K>0$ such that
$\parallel s.x\parallel \leq K.\parallel x
\parallel.\omega(s)$, for each $s\in S$, and $ x\in X$. But \[\parallel
s.x \parallel=\parallel \delta_s.x \parallel \leq \parallel \delta_s
\parallel_\omega
\parallel x \parallel=\omega(s).\parallel x \parallel,\] and the proof
is complete.\hfill$\blacksquare~$

We shall regard an essential Banach $M_a(S,\omega)$-module as a
unital Banach $M(S,\omega)$-module. Let $E$ be a unital Banach
$M_a(S,\omega)$-module. Then $\delta_s.x$ and $x.\delta_s$ are
defined in $E$ for each $s\in S$ and $x\in E$, they are often
denoted by $s.x$ and $x.s$ respectively. Similarly we define
$s\times\mu$ and $\mu\times s$ in $E^\prime$ for each $s\in S$ and
$\mu\in E^\prime$.\\ We need the following lemma for the next
results.
\paragraph{\bf Lemma 2.8.}Let $S$ be a locally compact foundation
semigroup with identity and  let $\omega$ be a weight function on
$S$.  we have \[\int_S\delta_{st}
 d\nu(s)=\nu*\delta_t  \qquad(\nu\in
 M_a(S,\omega), t\in S).\]
\paragraph{\bf Proof.} Let $t\in S$. For each $g\in C_0(S,\omega^{-1})$
we have,
\[<g,\int_S \delta_{st} d\nu(s)>=\int_S
g(st)d\nu(s)=<g,\nu*\delta_t>,\] and this gives the result.
\hfill$\blacksquare~$

 We establish a correspondence   between
continuous derivations from $M_a(S,\omega)$ and $S$-derivations into
dual modules.
\paragraph{\bf Lemma 2.9.}Let $E$ be a weighted $S$-module and every
weighted $S$-derivation from $S$ into $E^\prime$ be inner. Then
$H^1(M_a(S,\omega),E^\prime)=\{0\}$.
\paragraph{\bf Proof.}Let $D:M_a(S,\omega)\longrightarrow E^\prime$ be
a continuous derivation, where $E$ is a neo-unital Banach
$M_a(S,\omega)$-module. By 2.7, $E$ is a weighted $S$-module.
$M_a(S,\omega)$ has a bounded approximate identity, and then by 2.5,
there exists a unique extension derivation of $D$ on $M(S,\omega)$.
We denote this extension by $D$. Define $\varrho:S\longrightarrow
E^\prime$ by $\varrho(s)=D(\delta_s)$. We claim that $\varrho$ is a
weighted $S$-derivation. Let $s,t\in
S$,\[\varrho(st)=D(\delta_{st})=D(\delta_s*\delta_t)=D(\delta_s).\delta_t+\delta_s.D(\delta_t)=\varrho(s).t+s.\sigma(t).\]
But \[\frac{\parallel \varrho(s)
\parallel}{\omega(s)}=\frac{\parallel D(\delta_s)
\parallel}{\omega(s)}\leq \frac{\parallel D \parallel \parallel
\delta_s
\parallel_\omega}{\omega(s)}=\parallel D \parallel,\] and then
\[\sup\{\frac{\parallel \varrho(s)
\parallel}{\omega(s)}:s\in S\}<\infty.\] Therefore $\varrho$ is a
weighted $S$-derivation and there is $\varphi \in E^\prime$ such
that $\varrho(s)=s.\varphi-\varphi.s$. Now if $f\in
l^1(S,\omega)\subseteq M(S,\omega)$, then $f=\Sigma_s
\lambda_s.\delta_s$ where $\parallel f \parallel_\omega=\Sigma_s
\mid \lambda_s\mid.\omega(s)<\infty$, then
\begin{eqnarray*}
D(f)&=&D(\Sigma_s \lambda_s.\delta_s)=\Sigma_s
\lambda_s.D(\delta_s)\\
&=&\Sigma_s \lambda_s.\varrho(s)=\Sigma_s
\lambda_s.(s.\varphi-\varphi.s)\\
&=&(\Sigma_s \lambda_s.\delta_s).\varphi-\varphi(\Sigma_s
\lambda_s.\delta_s)\\
&=&f.\varphi-\varphi.f,
 \end{eqnarray*}
 and hence by 2.4, $D$ is
inner and proof is complete.\hfill$\blacksquare~$
\paragraph{\bf Lemma 2.10.} Let $S$ be locally compact foundation
semigroup with identity and let $\omega$ be a weight function on
$S$. Let $E$ be a neo-unital Banach $M_a(S,\omega)$-module, and
$\varrho:S\longrightarrow E^\prime$ be weighted $S$-derivation
defined in lemma 2.9, then the map $\varrho:S\longrightarrow
(E^\prime,\sigma)$ is continuous.
\paragraph{\bf Proof.} Let $s_\alpha\longrightarrow s$ in $S$. For each
$\nu\in M_a(S,\omega)$ the maps $s\longrightarrow \nu*\delta_s$ and
$s\longrightarrow \delta_s*\nu$ are continuous, so
$\nu*\delta_{s_\alpha}\longrightarrow \nu*\delta_s$ and then by
above definition $\delta_{s_\alpha}\longrightarrow \delta_s$ in (so)
topology. But $D$ is continuous
 in $w^*$-topology and hence $D(\delta_{s_\alpha})\longrightarrow
D(\delta_s)$ in $w^*$-topology. Therefore
$\varrho(s_\alpha)\longrightarrow \varrho(s)$ in
$(E^\prime,\sigma)$. \hfill$\blacksquare~$\\
 Let $E$ be a weighted
$S$-module and $\varrho:S\longrightarrow E^\prime$ be a weighted
$S$-derivation.
 By lemma 2.7, $E$ is a neo-unital Banach $M_a(S,\omega)$-module. Set
$\varrho(s)=D(\delta_s)\quad(s\in S).$  The following lemma shows
how $D$ can be recovered from $\varrho$.
\paragraph{\bf Lemma 2.11.} Let $D$ and $\varrho$ be as above. Then
\[<x,D\nu>=\int_S <x,\varrho(s)>d\nu(s)\quad(\nu\in
M_a(S,\omega),x\in E).\]
\paragraph{\bf Proof.}Since the function $s\longrightarrow
<x,\varrho(s)>$ is continuous and, for each $\nu\in M_a(S,\omega)$
\begin{eqnarray*}
\mid<
x,D\nu>\mid&=&\mid\int_S<x,\varrho(s)>d\nu(s)\mid\leq\int_S\mid<x,\varrho(s)>\mid
d\mid\nu\mid(s)\\
 &\leq& \parallel x \parallel
\int_S\parallel\varrho(s)\parallel d\mid\nu\mid(s)\leq
\parallel x \parallel.M.\int_S\omega(s)d\mid\nu\mid(s)\\
&=&M.\parallel x \parallel.\parallel \nu \parallel_\omega<\infty,
\end{eqnarray*}
 for some $M>0$, the integrals are defined on $S$.
 We first claim that $$\int_S <\delta_t
.y,\varrho(s)>d\nu(s)= <\delta_t .y,D\nu>\quad(\nu\in
M_a(S,\omega),t\in S,y\in E).$$ Let $\nu\in M_a(S,\omega),t\in S$
and let $y\in E$. Then by  2.8, we have
\begin{eqnarray*}\int_S
<\delta_t.y,\varrho(s)>d\nu(s)&=& \int_S
<\delta_t.y,D(\delta_s)>d\nu(s)= \int_S
<y,D(\delta_s).\delta_t>d\nu(s)\\
&=&
\int_S<y,D(\delta_{st})>d\nu(s)-\int_S<y.\delta_s,D(\delta_t)>d\nu(s)\\
&=&
<y,D(\int\delta_{st}d\nu(s)>-<y.\int_S\delta_sd\nu(s),D(\delta_t)>\\
&=& <y,D(\nu*\delta_t)>-<y.\nu,D(\delta_t)>\\
&=& <y,D\nu\times \delta_t>= <\delta_t.y,D\nu>.
\end{eqnarray*}
 Now take $x\in E$. Then $x=\mu.y\in E$, where
$\mu\in M_a(S,\omega)$ and $y\in E$. We have
\begin{eqnarray*}
\int_S<x,\varrho(s)>d\nu(s)&=&
\int_S<\int_S\delta_td\mu(t).y,\varrho(s)>d\nu(s)\\
&=& \int_S\int_S<\delta_t.y,\varrho(s)>d\nu(s)d\mu(t)\\
&=& \int_S<\delta_t.y,D\nu>d\mu(t)= <\mu.y,D\nu>,
\end{eqnarray*}
and so the lemma follows. \hfill$\blacksquare~$

Let $\varrho:S\longrightarrow E^*$ be a $S$-derivation, and let
$D:(M_a(s,\omega),\|.\|_1)\longrightarrow (E^\prime,\|.\|)$ defined
by
$$<x,D\nu>=\int_S <x,\varrho(s)>d\nu(s)\quad(\nu\in M_a(S,\omega), x\in
E).$$ Then $D$ is continuous. So we have the following lemma.
\paragraph{\bf Lemma 2.12.} The above
map $D$ is a derivation.
\paragraph{\bf Proof.}Take $\nu,\mu\in M_a(S,\omega)$ and $x\in E$.
Then
\begin{eqnarray*}
<x,D(\nu*\mu)>&=&\int_S<x,\varrho(s)>d(\nu*\mu)(s)=
\int_S\int_S<x,\varrho(st)>d\nu(s)d\mu(t)\\
&=& \int_S\int_S<x,\varrho(s).t+s.\varrho(t)>d\nu(s)d\mu(t)\\
&=&
\int_S\int_S(<\delta_t.x,\varrho(s)>+<x.\delta_s,\varrho(t)>)d\nu(s)d\mu(t).
\end{eqnarray*}
We also have
\begin{eqnarray*}
<x,D\nu \times\mu>&=&<\mu.x,D\nu>=
\int_S<\mu.x,\varrho(t)>d\nu(t)\\
&=&\int_S<\int_S\delta_sd\mu(s).x,\varrho(t)>d\nu(t)\\
&=& \int_S\int_S<\delta_s.x,\varrho(t)>d\mu(s)d\nu(t),
\end{eqnarray*}
by a similar calculation  for $<x,\nu \times D\mu>$, it follows that
\[<x,D(\nu*\mu)>=<x,D\nu
\times\mu>+<x,\nu \times D\mu>.\]Thus $D$ is a derivation.
\hfill$\blacksquare~$

From the $S$-derivation $\varrho$, we obtained the derivation $D$;
define $\varrho^\prime$ from $D$ as before. By lemma 2.11, for each
$\nu\in M_a(S,\omega)$ and $x\in E$, we have
\[\int_S<x,\varrho(s)>d\nu(s)=<x,D\nu>=\int_S<x,\varrho^\prime(s)>d\nu(s),\]
and so the functions $s\longmapsto <x,\varrho(s)>$ and $s\longmapsto
<x,\varrho^\prime(s)>$ are equal as elements of
$L^\infty(S;M_a(S,\omega))$. But both of the functions are
continuous on $S$, and so $<x,\varrho(s)>=<x,\varrho^\prime(s)>$ for
every $s\in S.$ Hence $\varrho=\varrho^\prime$. Suppose that
$\varrho:S\longrightarrow E^\prime$ is an inner $S$-derivation, and
take $\xi\in E^\prime$ with $\varrho(s)=s \times \xi-\xi \times s,$
 foe every $s\in S$. Let $D$ be defined from $\varrho$. Then, for
$\nu\in M_a(S,\omega)$ and $x\in E$, we have
\begin{eqnarray*}
<x,D\nu>&=& \int_S<x,\varrho(s)>d\nu(s)= \int_S<x,s \times \xi-\xi
\times s>d\nu(s)\\
&=& \int_S<x,s \times \xi>d\nu(s)-\int_S<x,\xi \times s>d\nu(s)\\
&=& \int_S<x.\delta_s,\xi>d\nu(s)-\int_S<\delta_s.x,\xi>d\nu(s)\\
&=& <x.\int_S\delta_sd\nu(s),\xi>-<\int_S\delta_sd\nu(s).x,\xi>\\
&=& <x.\nu,\xi>-<\nu.x,\xi>= <x,\nu\times\xi>-<x,\xi\times\nu>\\
&=& <x,\nu\times\xi-\xi\times\nu>,
\end{eqnarray*}
and so $D\nu=\nu \times \xi-\xi \times \nu$ and $D$ is inner. We
have established the following result.
\paragraph{\bf Theorem 2.13.} Let $S$ be a locally compact foundation
semigroup with identity, $\omega$ be a weight function on $S$, and
let $E$ be an neo-unital Banach $M_a(S,\omega)$-module. Then the map
$D\longrightarrow \varrho$, where we define $\varrho(s)=D(\delta_s)$
($s\in S$), establishes a bijective correspondence between
continuous derivations $D:M_a(S,\omega)\longrightarrow E^\prime$ and
$S$-derivations $\varrho:S\longrightarrow E^\prime$. Further, inner
derivation correspond to inner $S$-derivations, and
$H^1(M_a(S,\omega),E^\prime)=\{0\}$ if and only if every
$S$-derivation into $E^\prime$ is inner.
\section{Application to weighted group algebras}
 Throughout, $G$ will be a locally compact group(that is, clearly a
foundation semigroup with identity) and $\lambda$ a fixed left Haar
measure on $G$. For a weight function $\omega$ on $G$, let
$L^1(G,\omega)$ be the Banach space of all complex functions $f$
such that
\[\parallel f
\parallel_\omega=\int_G \mid f(t)\mid
\omega(t)d\lambda(t)<\infty\,\] where as usual we equate function
$\lambda$ almost everywhere. Under convolution product
\[(f*g)(x)=\int_G f(xy^{-1})g(y)d\lambda(y)\quad(f,g\in
L^1(G,\omega),\quad x\in G),\]$L^1(G,\omega)$ becomes a Banach
algebras. When $\omega(t)\geq 1$, the algebra $L^1(G,\omega)$ is
called a Beurling algebra. Let $M(G,\omega)$ be the Banach space of
all complex regular Borel measures $\mu$ on $G$ such that $\parallel
\mu \parallel_\omega=\int_G \omega(t)d\mid\lambda\mid(t)<\infty.$ If
$C_0(G,\omega)$ is the Banach space of all functions $f$ on $G$ such
that, $f/\omega\in C_0(G)$ and $ \parallel f
\parallel=\sup_{x\in G}\mid f(x)/\omega(x)\mid,$ then by pairing
\[<\mu,\psi>=\int_G \psi(x)d\mu(x)\qquad (\mu\in M(G,\omega),
\psi\in C_0(G,\omega)),\] we have
$(C_0(G,\omega))^\prime=M(G,\omega)$, and we can define the product
of $\mu,\nu\in M(G,\omega)$ by
\[\int_G \psi(x)d(\mu*\nu)(x)=\int_G\int_G\psi(xy)d\mu(x)d\nu(y)
\qquad (\psi\in C_0(g,\omega)).\] On $M(G,\omega)$ we consider one
topology other than the norm topology and the weak topology; it is
the strong operator (so) topology in which a net $(\mu_\alpha)$
tends to a measure $\mu$ (so) if $\mu_\alpha*f\longrightarrow \mu*f$
for every $f\in L^1(G,\omega)$, in $\parallel.\parallel$- topology
[4].\\ It is well known that the algebra $L^1(G)$ has a bounded
approximate identity $(f_i)$ consisting of continuous functions with
compact support, where each $f_i$ contains the identity element $e$
of $G$. Moreover, we can assume that there exists a compact set $K$
containing the support of all the $f_i$. If $g\in L^1(G,\omega)$ is
with compact support $C$,  we have $C\subseteq CK$, whence,
\[\parallel g-g*f_i\parallel_\omega\leq (\sup \omega(x))_{x\in
CK}\parallel g-g*f_i\parallel_1\] and likewise
\[\parallel g-f_i*g\parallel_\omega\leq (\sup \omega(x))_{x\in
KC}\parallel g-f_i*g\parallel_1\] and since the set of functions
with compact support is dense in $L^1(G,\omega)$, then we have the
following theorem(see[6]).
\paragraph{\bf Theorem 3.1.}
The algebra $L^1(G,\omega)$ possesses a bounded approximate identity
$(f_i)$ such that each $f_i$ is a continuous function with compact
support and there exists a compact set $S$ containing the support of
all the $f_i$.
\paragraph{\bf Remark 3.2.}Let $G$ be a locally compact group and $S$
be a subset of $G$. Set $\tilde{\omega}(g)=\omega(g^{-1})$ and
$\omega^*(g)=\omega(g).\tilde{\omega}(g)$. Then $\omega$ is
diagonally bounded on $S$, if $\omega^*$ is bounded on $S$. Also, it
is clear that, a symmetric weight is diagonally bounded on $S$ if
and only if it is bounded on $S$[4].
\paragraph{\bf Corollary 3.3.} Let $G$ be a locally compact group and
let $\omega$ be a weight function on $G$. Let $E$ be an neo-unital
Banach $L^1(G,\omega)$-bimodule, then
$H^1(L^1(G,\omega),E^\prime)=\{0\}$ if every crossed homomorphism
from $G$ into $E^\prime$ is principle, and only if where $\omega^*$
is bounded.
\paragraph{\bf Proof.}
Let $H^1(L^1(G,\omega),E^\prime)=\{0\}$, then each continuous
 derivation $D:L^1(G,\omega)\longrightarrow E^\prime$ is inner. Let
$D:L^1(G,\omega)\longrightarrow E^\prime$ be a continuous
derivation. By theorem 2.13, this $D$ induces a weighted
$G$-derivation $\varrho$, that is inner. Set
$\Psi(s)=\varrho(s).s^{-1}$\quad($s\in G$). For each $s,t\in G$, we
have
\[\Psi(st)=\varrho(st).(st)^{-1}=[\varrho(s).t+s.\varrho(t)].t^{-1}.s^{-1}
=\varrho(s).s^{-1}+s.\varrho(t).t^{-1}.s^{-1}=\Psi(s)+s.\Psi(t).s^{-1}.\]
Furthermore, since $\varrho$ is $\omega$-bounded, there exists
$M_0>0$, such that $\parallel \varrho(s)
\parallel \leq M_0.\omega(s)$ for each $s\in G$, and since $\omega^*$
is bounded, there exists $M_1>0$ such that $\omega^*(s)\leq M_1$ for
each $s\in G$. Hence
\[\parallel \Psi(s)
\parallel \leq
\parallel \varrho(s) \parallel.\parallel \delta_{s^{-1}}
\parallel_\omega
\leq M_0.\omega(s).\omega(s^{-1})=M_0.\omega^*(s)\leq M_0.M_1=M,\]
for each  $s\in G$. This means that $\Psi$ is a crossed
homomorphism. But $\varrho$ is inner and clearly $\Psi$ becomes
inner. Conversely, let $\Psi:G\longrightarrow E^\prime$ be a crossed
homomorphism, and for every $s\in G$, set $\varrho(s)=\Psi(s).s$.
For each $s,t\in G$, we have
\[\varrho(st)=\Psi(st).(st)=[\Psi(s)+s.\Psi(t).s^{-1}].s.t=\varrho(s).t+s.\varrho(t).\]
Now suppose that there exists $M>0$ such that $\parallel \Psi(s)
\parallel\leq M$ for each $s\in G$, so for arbitrary $s\in G$,
\[\parallel \varrho(s) \parallel \leq \parallel \Psi(s)
\parallel.\parallel \delta_s \parallel_\omega \leq M.\omega(s)\]
Then $\varrho$ is a weighted $G$-derivation. But it is clear that
 $\varrho$ is inner if $\Psi$ is principle, and proof is complete if we
apply 2.13. \hfill$\blacksquare~$
\paragraph{\bf Theorem 3.4.} Let $G$ be a locally compact amenable
group, and let $\omega$ be a weighted on $G$. If $\omega^*$ is
bounded on $G$(Indeed, $\omega$ is diagonally bounded on $G$) then
$L^1(G,\omega)$ is amenable.
\paragraph{\bf Proof.} Let $E$ be a weighted $G$-module and
$\varrho:G\longrightarrow E^\prime$ be a weighted $G$-derivation.
Let $K$ be the $w^*$-closed convex hull of the set
$\{\varrho(g).g^{-1}:g\in G\}$. Since $\omega^*$ is bounded, there
is $M_0>0$ such that $\omega^*(g)\leq M_0$ for all $g\in G$, and
since $\varrho$ is a weighted $G$-derivation, there is $M>0$ such
that $\parallel \varrho(g) \parallel\leq M.\omega(g)$ for all $g\in
G, x\in E$. Then for all $g\in G$,\[\parallel \varrho(g).g^{-1}
\parallel \leq
\parallel \varrho(g) \parallel.\parallel \delta_g^{-1} \parallel_\omega
\leq M.\omega(g).\omega(g^{-1}) =
 M.\omega^*(g) \leq M.M_0,\] and hence by Banach-Alaoglu theorem $K$ is
 compact. Now we define
$g\times \varphi=g.\varphi.g^{-1}+\varrho(g).g^{-1}$ where $(g\in G,
\varphi\in E^*)$. We have
\begin{eqnarray*}
gh\times\varphi&=&gh.\varphi.(gh)^{-1}+\varrho(gh).(gh)^{-1}\\
&=&gh.\varphi.h^{-1}.g^{-1}+\varrho(g).g^{-1}+g.\varrho(h).h^{-1}.g^{-1}\\
&=&g.(h.\varphi.h^{-1}+\varrho(h).h^{-1}).g^{-1}+\varrho(g).g^{-1}\\
&=&g\times(h\times\varphi),
\end{eqnarray*}
and further $\varrho(e)=0$. Then $\times$ is a group action. It is
clear that $G$ acts affinely on $E^*$ by this action. Let
$(g_\alpha)$ be a net in $G$ such That $g_\alpha\longrightarrow g$
in $G$. Since $g\mapsto g.x$ is continuous and, there is $M>0$ such
that $\parallel g.x
\parallel\leq M.\parallel x
\parallel.\omega(g)$ for all $g\in G, x\in E$ and $\omega$ is
continuous , we have
\[<x,g_\alpha.\varphi.g_\alpha^{-1}>=<g_\alpha^{-1}.x.g_\alpha,\varphi>\longrightarrow
<g^{-1}.x.g,\varphi>=<x,g.\varphi.g^{-1}>\] and
\[\mid<g_\alpha^{-1}.x,\varrho(g_\alpha)>-<g^{-1}.x,\varrho(g)>\mid
\leq
\parallel \varrho(g)-\varrho(g_\alpha)\parallel.\parallel
g_\alpha^{-1}.x\parallel +\parallel \varrho(g)\parallel.\parallel
g_\alpha^{-1}.x-g^{-1}.x\parallel\]\[\leq M.\parallel x
\parallel.\omega(g_\alpha^{-1}).\parallel
\varrho(g)-\varrho(g_\alpha)\parallel+
\parallel \varrho(g)\parallel.\parallel
g_\alpha^{-1}.x-g^{-1}.x\parallel\longrightarrow 0.\] and hence
$g_\alpha \times \varphi\longrightarrow g \times \varphi$ in
$w^*$-topology. Now let $\varphi_\alpha\longrightarrow \varphi$ in
$w^*$-topology. Then
\[<x,g.\varphi_\alpha.g^{-1}>=<g^{-1}.x.g,\varphi_\alpha>\longrightarrow
<g^{-1}.x.g,\varphi>=<x,g.\varphi.g^{-1}>,\] i.e $g \times
\varphi_\alpha\longrightarrow g\times\varphi$. This means that the
group action $\times$ is continuous in the first variable with
respect to the given topology on $G$ and continuous in the second
variable with respect to $w^*$-topology on $E'$. Then by Day's fixed
point theorem there is $\varphi \in K$ such that $g \times
\varphi=\varphi$ and hence
\[g.\varphi.g^{-1}+\varrho(g).g^{-1}=\varphi\Longrightarrow
D(g)=\varphi.g-g.\varphi.\] Thus by corollary 2.13, $L^1(G,\omega)$
is amenable.\hfill$\blacksquare~$\\

\end{document}